\documentclass[11pt]{article}
\usepackage{amssymb,amsfonts,amsmath,latexsym,epsf,tikz,url}
\newtheorem{theorem}{Theorem}[section]

\newtheorem{corollary}[theorem]{Corollary}
\newtheorem{lemma}[theorem]{Lemma}
\newtheorem{rem}[theorem]{Remark}

\newcommand{\proof}{\noindent{\bf Proof.\ }}
\newcommand{\qed}{\hfill $\square$\medskip}

\begin{document}

\title{On the anti-forcing number of graph powers}

\author{Neda Soltani$^{}$\footnote{Corresponding author} \and Saeid Alikhani}

\date{\today}

\maketitle

\begin{center}

Department of Mathematics, Yazd University, 89195-741, Yazd, Iran\\
{\tt Neda\_soltani@ymail.com}\\

\end{center}

\begin{abstract}
Let $G=(V,E)$ be a simple connected graph. A perfect matching (or Kekul\'e structure in chemical literature) of $G$ is a 
set of disjoint edges which covers all vertices of $G$. The anti-forcing number of $G$ is the smallest number of edges such that the remaining graph obtained by deleting these edges has a unique perfect matching and is denoted by $af(G)$. 
For every $m\in\mathbb{N}$, the $m$th power of $G$, denoted by $G^m$, is a graph with the same vertex set as $G$ such that two vertices are adjacent in $G^m$ if and only if their distance is at most $m$ in $G$. In this paper, we study 
the anti-forcing number of the powers of some graphs.

\end{abstract}


\medskip
\noindent{\bf AMS Subj.\ Class.:} 05C70, 05C76.

\section{Introduction}
Let $G$ be a simple graph with vertex set $V(G)$ and edge set $E(G)$. 
A matching $M$ in a graph $G$ is a collection of edges of $G$ such that no two edges from $M$ share a vertex.
If each vertex of $G$ is incident with exactly one edge in $M$, the matching $M$ is called perfect. In other words, all vertices in a perfect matching are saturated. Obviously, only graphs on an even number of vertices can have perfect matching. The study of perfect matchings 
has a long history in both mathematical and chemical literature. For more details on perfect matching, we refer the 
reader to see \cite{Lov}.

If $G$ is a graph that admits a perfect matching $M$, then a forcing set of $M$ is a subset $S$ of $M$ contained in no other perfect matchings of $G$. The minimum possible cardinality of forcing set $S$ is called the forcing number of $M$ and denoted by $f(G)$.

The notions of a forcing edge and the forcing number of a perfect matching first appeared in 1991 by Harary, Klein and \^Zivkovi\'c (\cite{Har}). The root of these concepts can be traced to the works \cite{{Kle}, {Ran}} by Randi\'c and Klein where the forcing number was introduced under the name of “innate degree of freedom” of a
Kekulé structure, which plays an important role in the resonance theory in chemistry.

Over the past two decades, more and more mathematicians were attracted to the study on forcing sets and the forcing numbers of perfect matchings of a graph. The scope of graphs in consideration has been extended from polyhexes to various bipartite graphs and non-bipartite graphs.
Zhang, Ye and Shiu proved that the forcing number of fullerenes has a lower bound three and there are infinitely many fullerenes achieving this bound (\cite{Zha}). 

In 2007, Vukicevi\'c and Trinajsti\'c in \cite{{Vu},{Tr}} introduced the anti-forcing number of a graph $G$. For $S\subseteq E(G)$, let $G-S$ denote the graph obtained by removing $S$ from $G$. Then $S$ is called an anti-forcing set if $G-S$ has a unique perfect matching. The cardinality of the smallest anti-forcing set is called the anti-forcing number of $G$ and denoted by $af(G)$. An edge $e$ of $G$ is called an anti-forcing edge if $G-e$ has a unique perfect matching. Note that $af(G)=|E(G)|$ if and only if $G$ does not have any perfect matching. A graph $G$ is called odd and even graph, if the number of vertices of $G$ is odd and even, respectively. 

Yang, Zhang and Lin proved that the anti-forcing number of every fullerene has a lower bound four. They also showed that for every even $n\geq 20$ $(n\neq 22, 26)$, there exists a fullerene with $n$ vertices that has the anti-forcing number four, and the fullerene with $26$ vertices has the anti-forcing number five (\cite{Yan}).

Recently, Lei, Yeh and Zhang in \cite{Lei} defined the anti-forcing number of a perfect matching $M$ of a graph $G$ as the minimal number of edges not in $M$ whose removal to make $M$ as a single perfect matching of the resulting graph and denoted by $af(G,M)$. By this definition, the anti-forcing number of a graph $G$ is the smallest anti-forcing number over 
all perfect matchings of $G$. They also proved that the anti-forcing number has a close relation with the
forcing number. In fact, for any perfect matching $M$ of a graph $G$, $f(G, M)\leq af(G, M)\leq (\Delta-1)f(G, M)$ where 
$\Delta$ denotes the maximum degree of the graph $G$.

For every positive integer $k$, the $k$-power of $G$ is defined on the $V(G)$ by adding edges joining any two distinct vertices $x$ and $y$ with distance at most $k$ in $G$ and is denoted by $G^k$ (\cite{{Agn},{Kra}}). In other words, $E(G^k)=\{xy: 1\leq d_G(x, y)\leq k\}$. The following lemma follows from the definition of the power of graphs.

\begin{lemma}\cite{Dis}\label{int}
Let $G$ be a connected graph of order $n$ and diameter $d$. 
\begin{enumerate}
\item[(i)] For every natural number $t \geq d$, $G^t=K_n$.
\item[(ii)] (Theorem 1 in \cite{Hob}) Let $k=mn$ where $m, n\in \Bbb{N}$. Then $G^k={(G^m)}^n$.
\item[(iii)] (Lemma 2.1 in \cite{An}) Let $x$ and $y$ be two vertices of $G$. Then 
$d_{G^k}(x, y)=\left\lceil\frac{d_G(x, y)}{k}\right\rceil$.
\end{enumerate}
\end{lemma}

The anti-forcing number of some specific graphs such as paths, cycles, friendships and cactus graphs has studied in \cite{anti}. In this paper, we consider the anti-forcing number of the powers of these graphs. 

\section{Anti-forcing  number of power of specific graphs}

In this section, we compute the anti-forcing number of the powers of some certain graphs such as paths, cycles, friendship graphs and also some graphs that are of importance in Chemistry. The following theorem gives the anti-forcing number of the powers of even paths.

\begin{theorem}\label{afp}
Let $P_k$ be a path of order even $k$ and $m\geq 2$. 

\begin{enumerate}
\item[(i)] If $k=2m$ \textbf{or} $k<2m$ and $d(x_{k-m+1}, x_m)<2$, then $af(P_k^m)=(k-m)(m-1)$.
\item[(ii)] If $k<2m$ and $d(x_{k-m+1}, x_m)\geq 2$, then $$af(P_k^m)=(k-m)(m-1)+\sum_{i=1}^{\left\lfloor{\frac{2m-k-1}{2}}\right\rfloor}(2m-k-2i).$$
\item[(iii)] If $k>2m$, then $$af(P_k^m)=\left\lfloor{\frac{k}{2m}}\right\rfloor m(m-1)+af(P_{k-2m\left\lfloor{\frac{k}{2m}}\right\rfloor}^m).$$
\end{enumerate}
\end{theorem}

\proof
Suppose that $\{v_1, v_2, ..., v_k\}$ is the vertex set of $P_k^m$ and the set $S$ has the smallest cardinality over all anti-forcing sets of the graph $P_k^m$. For every $1\leq i\leq k-m$, there exist $m$ possibilities to saturate the vertex $v_i$ and so to have an anti-forcing set with the smallest cardinality, we put the edges $v_iv_{i+1}, v_iv_{i+2}, ..., v_iv_{m+i-1}$ in $S$. Thus $|S|=\sum_{i=1}^{k-m}(m-1)$. Note that for every $1\leq i\leq m$, the edges $v_iv_{m+i}\in M$ where $M$ is the unique perfect matching of $P_k^m-S$. 
If $k=2m$, then $$af(P_k^m)=|S|=\sum_{i=1}^{m}(m-1).$$ Otherwise, assume that $k<2m$. Then the vertices $v_{k-m+1}, v_{k-m+2}, ..., v_{m-2}, v_{m-1}, v_m$ are $M$-unsaturated. In this case, to have the smallest anti-forcing set, we add the edges $v_{k-m+i}v_{k-m+i+1}, v_{k-m+i}v_{k-m+i+2}, ..., v_{k-m+i}v_{m-i}$ to $S$ where $1\leq i\leq \left\lfloor{\frac{2m-k-1}{2}}\right\rfloor$. It is easy to see that the number of these edges is $2m-k-2i$ and so $$af(P_k^m)=|S|=\sum_{i=1}^{k-m}(m-1)+\sum_{i=1}^{\left\lfloor{\frac{2m-k-1}{2}}\right\rfloor}(2m-k-2i).$$ Note that if there just exist two $M$-unsaturated vertices, then the second summation is zero. 

Now suppose that $k>2m$. Then the vertices $v_{2m+1}, v_{2m+2}, ..., v_k$ are $M$-unsaturated. In this case, we should divide the graph $P_k^m$ into subgraphs of order $2m$. Then we have $\left\lfloor{\frac{k}{2m}}\right\rfloor$ subgraphs $P_{2m}^m$ and a subgraph with $k-2m\left\lfloor{\frac{k}{2m}}\right\rfloor$ vertices. Since $k-2m\left\lfloor{\frac{k}{2m}}\right\rfloor<2m$, so by the previous cases, we have the result. 
\qed

\begin{rem}
Clearly, the paths of order odd $k$ does not have any perfect matching and so $af(P_k^m)=|E(P_k^m)|$. It is easy to 
see that for every $m<k-1$, there exist $k-m$ pairs of vertices which their distance in $P_k$ is $m$. So we have 
$$af(P_k^m)=|E(P_k^m)|=|\{xy: 1\leq d_G(x, y)\leq m\}|=\sum_{i=1}^{m}(k-i)=mk-\frac{m}{2}(m+1).$$ Otherwise, $P_k^m$ is a complete graph of order $k$ and so $af(P_k^m)=\frac{k(k-1)}{2}$.
\end{rem}

\medskip
We need the following easy lemma to obtain a result for the anti-forcing number of the powers of even cycles. 

\begin{lemma}\label{power}
Let $G$ be a connected graph of order $n$. For every $m_1, m_2 \in \Bbb{N}$, if $m_1<m_2$, then $af(G^{m_1})\leq af(G^{m_2})$.
\end{lemma}

\begin{theorem}\label{afc}
Let $C_k$ be a cycle of order even $k$. Then for every $m\geq 2$,
$$\frac{k+8}{4}\leq af(C_k^m)\leq \frac{k(k-2)}{4}.$$
\end{theorem}

\proof
By Lemma \ref{power}, for every $m\geq 2$, $af(C_k^m)\geq af(C_k^2)$.
Suppose that $V(C_k)=\{v_1, v_2, v_3, ..., v_k\}$ and $S_1$ has the smallest cardinality over all anti-forcing sets of 
graph $C_k^2$. So $af(C_k^2)=|S_1|$.
Clearly, we have four possibilities to saturate the vertex $v_1$ and to have an anti-forcing set with the smallest cardinality, we shall put three edges $v_1v_2, v_1v_3$ and $v_1v_{k-1}$ into $S_1$ and so $v_1v_k\in M_1$ where $M_1$ is the unique perfect matching of $C_k^2-S_1$. Then there exist two edges $v_2v_3$ and $v_2v_4$ to saturate the vertex $v_2$. We should add the edge $v_2v_3$ into $S_1$ and so $v_2v_4, v_3v_5\in M_1$. As the same way, we can put the edge $v_6v_7$ into $S_1$ and so $v_6v_8, v_7v_9\in M_1$. By continuing this method, we have
$$S_1=\{v_1v_2, v_1v_3, v_1v_{k-1}, v_2v_3, v_6v_7, v_{10}v_{11}, v_{14}v_{15}, ..., v_{k-4}v_{k-3}\}.$$
So $|S_1|=3+l$ and $|M_1|=1+2l$ where $l$ is the number of edges $v_2v_3, v_6v_7, ...,v_{k-4}v_{k-3}$. In addition, for every even $k$, the number of edges in a perfect matching of a graph of order $k$ is equal to $\frac{k}{2}$. If 
$\frac{k}{2}$ is odd, then $\frac{k}{2}=1+2l$ and so $l=\frac{k-2}{4}$. Thus we have $|S_1|=3+l=3+\frac{k-2}{4}= \frac{k+10}{4}$. Otherwise, there exists an edge $v_{k-2}v_{k-1}$ that we shall add into $M_1$. Thus
$|M_1|=2+2l$ and so $l=\frac{k-4}{4}$. Therefore $|S_1|=3+\frac{k-4}{4}=\frac{k+8}{4}$ and so for every $m\geq 2$,
$$af(C_k^m)\geq af(C_k^2)\geq \frac{k+8}{4}.$$ 

Now assume that $S_2$ has the smallest cardinality over all anti-forcing sets of graph $C_k^{\frac{k}{2}}$. Thus $af(C_k^{\frac{k}{2}})=|S_2|$. Since the diameter of $C_k$ is equal to $\frac{k}{2}$, so for every $m\geq \frac{k}{2}$, $C_k^m$ is a complete graph of order $k$. Thus the upper bound of $af(C_k^m)$ is equal to $af(C_k^{\frac{k}{2}})$.
Clearly, we have $k-1$ possibilities to saturate the vertex $v_1$. Without loss of generality, we assume that the edge $v_1v_k\in M_2$ where $M_2$ is the unique perfect matching of $C_k^{\frac{k}{2}}-S_2$.
So we can put $k-2$ other edges into $S_2$. The vertex $v_2$ can be saturated by one of the edges $v_2v_3, v_2v_4, ..., v_2v_{k-2}$ or $v_2v_{k-1}$. Suppose that $v_2v_{k-1}\in M_2$. Then $k-4$ other edges belong to $S_2$. As the same way, for every $1\leq i\leq \frac{k}{2}-1$, there exist $k-2i+1$ possibilities to saturate the vertex $v_i$ and to have an anti-forcing set with the smallest cardinality, we shall put the edges $v_iv_{i+1}, v_iv_{i+2}, ..., v_iv_{k-i}$ into $S_2$. Therefore for every $m\geq 2$, 
\begin{equation*}
af(C_k^m)\leq af(C_k^{\frac{k}{2}})=|S_2|
=\sum_{i=1}^{\frac{k}{2}-1}(k-2i)=\frac{k(k-2)}{4}.~~~~\square
\end{equation*} 
\begin{rem}
\begin{enumerate}
\item[(i)] The lower and uppor bounds of Theorem \ref{afc} are sharp. 
\item[(ii)]  Clearly, the cycles of order odd $k$ does not have any perfect matching and thus $af(C_k^m)=|E(C_k^m)|$. 
If $m<\left\lfloor\frac{k}{2}\right\rfloor$, then there exist $k$ pairs of vertices which their distance in $C_k$ is $m$ and so $af(C_k^m)=mk$. Otherwise, $C_k^m$ is a complete graph and so $af(C_k^m)=k\left\lfloor\frac{k}{2}\right\rfloor$ where $\left\lfloor\frac{k}{2}\right\rfloor$ is the diameter of $C_k$.
\end{enumerate}
\end{rem}

\begin{figure}
	\begin{center}
	\includegraphics[width=10cm,height=2.3cm]{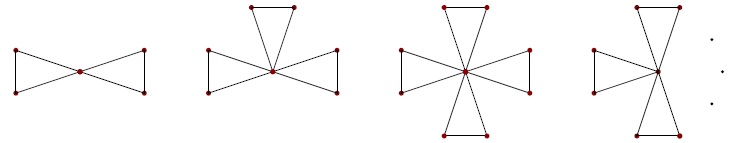}
	\caption{\small Friendship graphs $F_2, F_3, F_4$ and $F_k$, respectively.}
	\label{F_k}
\end{center}
\end{figure}

Here we consider the anti-forcing number of the powers of friendship graphs. The friendship (or Dutch-Windmill) graph 
$F_k$ is a graph can be constructed by joining $k$ copies of the cycle graph $C_3$ with a common vertex. Some examples of friendship graphs are shown in Figure \ref{F_k}. Clearly, a friendship graph $F_k$ has $2k+1$ vertices and $3k$ edges. 
Thus for every $m\in \Bbb{N}$, $af(F_k^m)=|E(F_k^m)|$ and since the diameter of $F_k^m$ is equal to 2, so for every 
$m\geq 2$, $F_k^m$ is a complete graph of order $2k+1$. Hence $$af(F_k^m)=2k^2+k.$$

\begin{figure}
	\begin{center}
	\includegraphics[width=6cm,height=1.5cm]{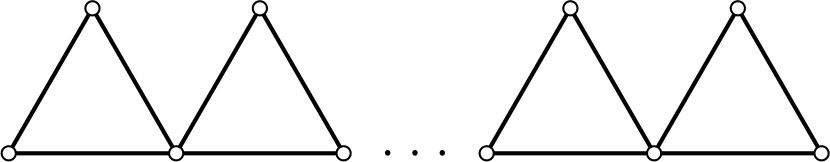}
	\caption{\small A chain triangular cactus $T_k$ .}
	\label{T_k}
\end{center}
\end{figure}

\medskip
Here we shall study the anti-forcing number of the powers of cactus graphs that are of importance in chemistry. 
A cactus graph is a connected graph in which no edge lies in more than one cycle. Consequently, each block of a cactus graph is either an edge or a cycle. If all blocks of a cactus $G$ are cycles of the same size $i$, the cactus is $i$-uniform. 
A triangular cactus is a graph whose blocks are triangles, i.e., a $3$-uniform cactus. 
A vertex shared by two or more triangles is called a cut-vertex. If each triangle of a triangular cactus $G$ has at most two cut-vertices and each cut-vertex is shared by exactly two triangles, we say that $G$ is a chain triangular cactus. The number of triangles in $G$ is called the length of the chain. Obviously, all chain triangular cactuses of the same length are isomorphic. Hence we denote the chain triangular cactus of length $k$ by $T_k$ (\cite{Jah}). An example of a chain triangular cactus is shown in Figure \ref{T_k}. Clearly, a chain triangular cactus of length $k$ has $2k+1$ vertices and 
$3k$ edges. So for every $m\geq 2$,

\medskip
\medskip
$\begin{array}{l}
\medskip
 
af(T_k^m)=|E(T_k^m)|=3k+4(k-1)+4(k-2)+...+4\big(k-(m-1)\big)\\
\medskip
~~~~~~~~~~=3k+4k(m-1)-2m(m-1)\\
\medskip
~~~~~~~~~~=4km-k-2m^2+2m.
\end{array}$

\begin{figure}
	\begin{center}
		\begin{minipage}{5cm}
			\hspace{5cm}
			\includegraphics[width=5cm,height=2cm]{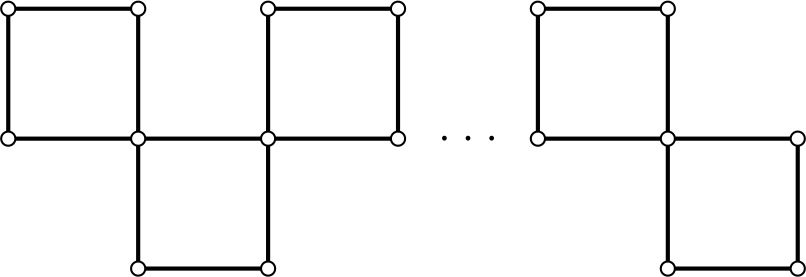}
		\end{minipage}
		\begin{minipage}{5cm}
			\includegraphics[width=6cm,height=1.7cm]{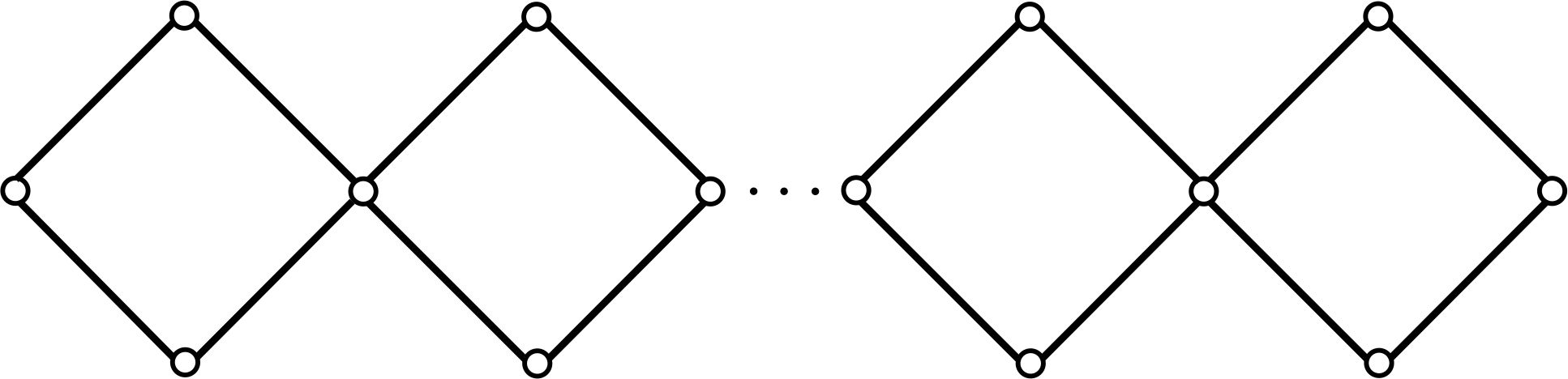}
		\end{minipage}
		\caption{ \label{1} \small An ortho-chain square cactus $O_k$ and a para-chain square cactus $Q_k$, respectively.}
	\end{center}
\end{figure}

\medskip
By replacing triangles in the definition of triangular cactus by cycles of length $4$, we obtain a cactus whose every block is $C_4$. We call such a cactus, square cactus. Note that the internal squares may differ in the way they connect to their neighbors. If their cut-vertices are adjacent, we say that such a square is an ortho-square and if the cut-vertices are not adjacent, we call the square a para-square (\cite{Jah}). An ortho-chain square cactus of length $k$, $O_k$, and a para-chain square cactus of length $k$, $Q_k$, are shown in Figure \ref{1}. 
The following theorem gives the anti-forcing number of the powers of the ortho-chain square cactuses of length even $k$.

\begin{theorem}\label{ortho}
Let $O_k$ be an ortho-chain square cactus of length $k$. If $k$ is even, then for every $m\geq 4$, 
$$af(O_k^m)=af(O_k^{m-1})+9(k-m)+15$$ and $af(O_k^2)=10k-4$, $af(O_k^3)=18k-16$. 
\end{theorem}

\begin{figure}
	\begin{center}
	\includegraphics[width=9cm,height=2.5cm]{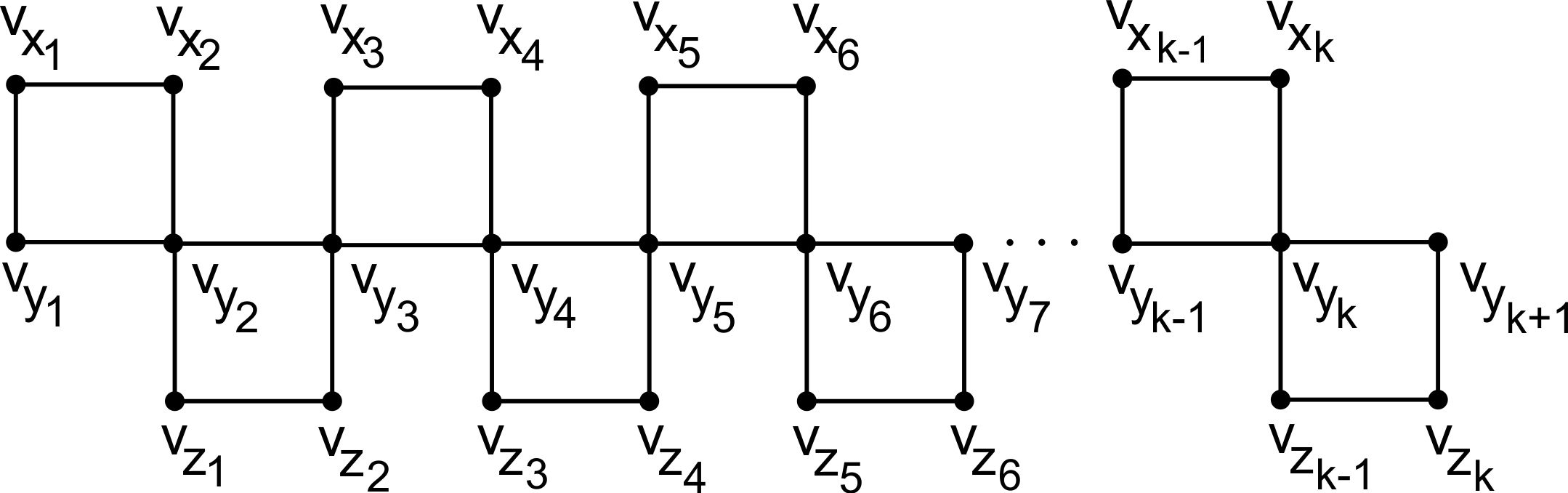}
	\caption{\small A labeled ortho-chain square cactus.}
	\label{O_k^m}
\end{center}
\end{figure}

\proof
Since $k$ is even, so the order of $O_k$ is odd. Thus we have 

\medskip
$\begin{array}{l}
\medskip
 
af(O_k^m)=|E(O_k^m)|=|\{xy:1\leq d_{O_k}(x, y)\leq m\}|\\
\medskip
~~~~~~~~~~=|\{xy:1\leq d_{O_k}(x, y)\leq m-1\}|+|\{xy:d_{O_k}(x, y)=m\}|\\
\medskip
~~~~~~~~~~=af(O_k^{m-1})+|\{xy:d_{O_k}(x, y)=m\}|
\end{array}$

\medskip
For particular cases $m=2, 3$, the proof is easy and left to the reader. Suppose that $m\geq 4$.
We shall compute the number of pairs of vertices which their distance is $m$. We label the vertices of $O_k^m$ with $v_{x_1}, v_{x_2}, ..., v_{x_k}, v_{y_1}, v_{y_2}, ..., v_{y_{k+1}}$ and $v_{z_1}, v_{z_2}, ..., v_{z_k}$ as shown in Figure \ref{O_k^m}. We have the following cases:

\medskip 
1) $d(v_{x_i}, v_{x_{m+i-2}})=m$, ~~~~~~~~~~~~~~~~~~~2) $d(v_{x_i}, v_{y_{m+i-1}})=m$,

3) $d(v_{x_i}, v_{z_{m+i-3}})=m$,~~~~~~~~~~~~~~~~~~~~4) $d(v_{y_i}, v_{x_{m+i-1}})=m$,

5) $d(v_{y_i}, v_{y_{m+i}})=m$,~~~~~~~~~~~~~~~~~~~~~~~6) $d(v_{y_i}, v_{z_{m+i-2}})=m$,

7) $d(v_{z_i}, v_{x_{m+i-1}})=m$,~~~~~~~~~~~~~~~~~~~~8) $d(v_{z_i}, v_{y_{m+i}})=m$,

9) $d(v_{z_i}, v_{z_{m+i-2}})=m$.

\medskip
\medskip
Here we compute the number of pairs of vertices which satisfy case 1. 
Since the number of vertices with index $x$ is $k$, so we have $i+m-2\leq k$ and thus $i$ can choose every element between $1$ and $k-m+2$. So there exist $k-m+2$ pairs in case 1. With similar approach, the number of pairs of vertices which satisfy another cases is $k-m+2, k-m+3, k-m+1, k-m+1, k-m+2, k-m+1, k-m+1$ and $k-m+2$, respectively and the result is obtained by adding these amounts and direct substitution into above formula.
\qed

\begin{corollary}\label{coro}
Let $O_k$ be an ortho-chain square cactus of length even $k$. Then for every $m\geq 4$, 
$$af(O_k^m)=9k(m-1)-\big(\frac{m-3}{2}\big)(9m+6)-16.$$
\end{corollary}

\proof
Let $X_m:=af(O_k^m)$. Then by Theorem \ref{ortho}, for every $m\geq 4$, we have
$X_m=X_{m-1}+9(k-m)+15$. Now by solving this recurrence relation, we have 
$$X_m=X_3+\big(\frac{m-3}{2}\big)(18k-9m-6)$$ and $X_3=18k-16$. So the result follows.
\qed

In the next theorem, we present the anti-forcing number of the powers of the para-chain square cactuses of length even 
$k$.

\begin{theorem}\label{para}
Let $Q_k$ be a para-chain square cactus of length $k$. If $k$ is even, then
    \begin{enumerate}
           \item [(i)] for every even $m\geq 4$, $$af(Q_k^m)=af(Q_k^{m-1})+5(k-\frac{m}{2})+1$$ and $af(Q_k^2)=10k-4$.
           \item [(ii)] for every odd $m\geq 3$, $$af(Q_k^m)=af(Q_k^{m-1})+4(k-\left\lfloor\frac{m}{2}\right\rfloor).$$
    \end{enumerate}
\end{theorem}

\begin{figure}
	\begin{center}
	\includegraphics[width=7.5cm,height=2.3cm]{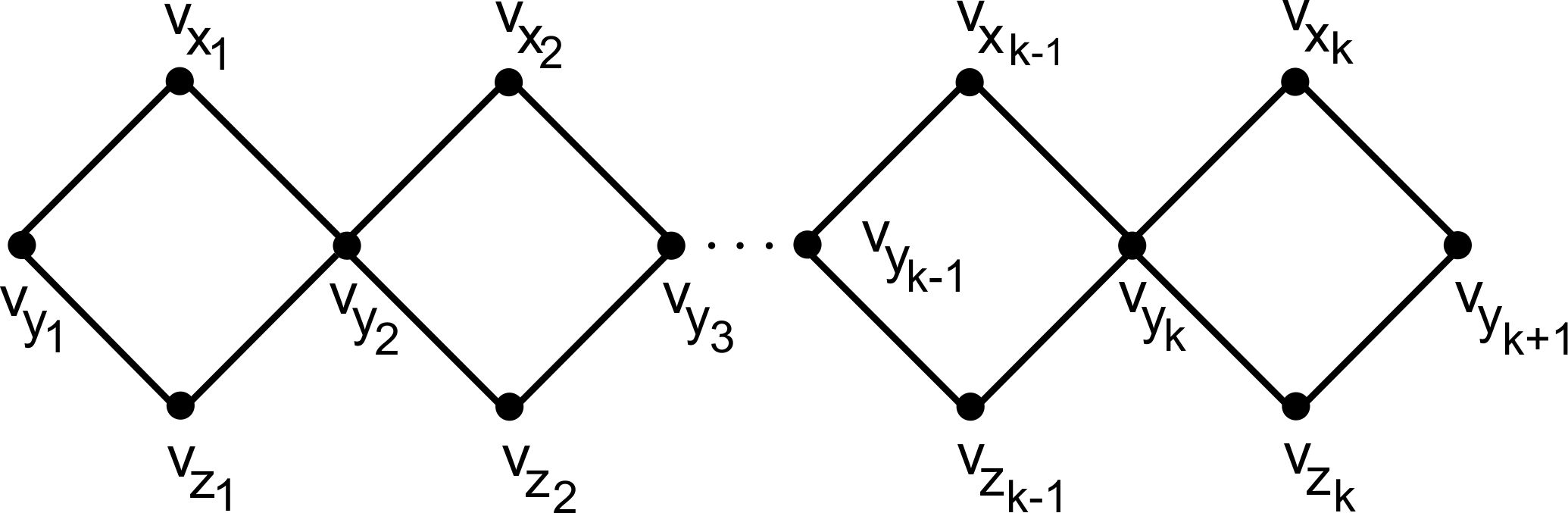}
	\caption{\small A labeled para-chain square cactus.}
	\label{Q_k^m}
\end{center}
\end{figure}

\proof
Since $k$ is even, so the order of $Q_k$ is odd. Thus in the same way as proof of  Theprem \ref{ortho}, it is sufficient to compute the number of pairs of vertices which their distance is $m$. We label the vertices of $Q_k^m$ with $v_{x_1}, v_{x_2}, ..., v_{x_k}, v_{y_1}, v_{y_2}, ..., v_{y_{k+1}}$ and $v_{z_1}, v_{z_2}, ..., v_{z_k}$ as shown in Figure \ref{Q_k^m}. Then
 \begin{enumerate}
           \item [(i)] for every even $m\geq 4$, we have the following cases:

1) $d(v_{x_i}, v_{x_{\frac{m}{2}+i}})=m$, ~~~~~~~~~~~~~~~~~~~2) $d(v_{x_i}, v_{z_{\frac{m}{2}+i}})=m$,

3) $d(v_{y_i}, v_{y_{\frac{m}{2}+i}})=m$, ~~~~~~~~~~~~~~~~~~~4) $d(v_{z_i}, v_{x_{\frac{m}{2}+i}})=m$,

5) $d(v_{z_i}, v_{z_{\frac{m}{2}+i}})=m$.

\medskip
The number of pairs of vertices which satisfy case 3 is $k-\frac{m}{2}+1$ and the number of pairs of vertices which satisfy another cases is the same and equal to $k-\frac{m}{2}$. So by adding these amounts, we have
$$af(Q_k^m)=af(Q_k^{m-1})+5(k-\frac{m}{2})+1.$$

Clearly, if $m=2$, then $af(Q_k^2)=10k-4$.

\item[(ii)] for every odd $m\geq 3$, we have the following cases:

1) $d(v_{x_i}, v_{y_{\left\lfloor\frac{m}{2}\right\rfloor+i+1}})=m$, ~~~~~~~~~~~~~~~~~~~2) $d(v_{y_i}, v_{x_{\left\lfloor\frac{m}{2}\right\rfloor+i}})=m$,

3) $d(v_{y_i}, v_{z_{\left\lfloor\frac{m}{2}\right\rfloor+i}})=m$, ~~~~~~~~~~~~~~~~~~~~~4) $d(v_{z_i}, v_{y_{\left\lfloor\frac{m}{2}\right\rfloor+i+1}})=m$.

\medskip
There exist $k-\left\lfloor\frac{m}{2}\right\rfloor$ pairs of vertices which satisfy each of these cases and we have the result with similar approach.
\qed
   \end{enumerate}

With similar method that used in the proof of Corollary \ref{coro}, we have the following corollary which gives the anti-forcing number of the powers of the para-chain square cactuses of length even $k$, directly.

\begin{corollary}
Let $Q_k$ be a para-chain square cactus of length even $k$.  
   \begin{enumerate}
           \item [(i)] For every even $m\geq 4$, 
$$af(Q_k^m)=\big(\frac{9m+2}{2}\big)k-\big(\frac{m-2}{8}\big)(9m+16)-4.$$
           \item [(ii)] For every odd $m\geq 3$, 
$$af(Q_k^m)=\big(\frac{9m+1}{2}\big)k-\big(\frac{m-1}{8}\big)(9m-11)-4.$$
   \end{enumerate}
\end{corollary}

\section{Conclusion} 
In this paper, we have considered the anti-forcing number of the powers of paths, cycles, friendships and chain triangular cactuses. In particular, we have established some formulas for the anti-forcing number of the powers of the chain square cactuses of length even $k$. It would be interesting to obtain the anti-forcing number of the powers of these graphs when their length is odd. Therefore we end this section by the following problem. \medskip\\
\textbf{Problem}: Find the anti-forcing number of the powers of the chain square cactuses of length $k$ for when $k$ is odd.



\end{document}